\documentclass[12pt]{amsart}

\usepackage{amsmath,amsfonts,amsthm}

\renewcommand{\theequation}{\thesection.\arabic{equation}}
\newtheorem{theorem}{Theorem}

\newcommand{\eqnsection}{
\renewcommand{\theequation}{\thesection.\arabic{equation}}
    \makeatletter
    \csname  @addtoreset\endcsname{equation}{section}
    \makeatother}
\eqnsection

\def\R{{\mathbb R}}

\def\E{{\mathcal E}}
\def\A{{\mathcal A}}
\def\L{{\mathcal L}}

\author[N. Enriquez]{Nathana\"el ENRIQUEZ}
\address{Laboratoire Modal'X, Universit\'e Paris 10, 200
Avenue de la R\'epublique, 92000 Nanterre, France}
\email{nenriquez@u-paris10.fr}

\author[C. Sabot]{Christophe SABOT}
\address{Universit\'e de
Lyon, Universit\'e Lyon 1, Institut Camille Jordan, CNRS UMR 5208,
43, Boulevard du 11 novembre 1918, 69622 Villeurbanne Cedex,
France} \email{sabot@math.univ-lyon1.fr}

\author[M. Yor]{Marc YOR}
\address{Laboratoire de Probabilit\'es et Mod\`eles Al\'eatoires, CNRS UMR 7599, Universit\'e Paris 6, 4
place Jussieu, 75252 Paris Cedex 05, France}
\email{deaproba@proba.jussieu.fr}


\title[Square-root boundaries]{Renewal series and  square-root boundaries for Bessel processes }

\begin{document}

\maketitle

{\footnotesize \noindent{\slshape\bfseries Abstract.} }
We show how a description of   Brownian exponential functionals  as a renewal series gives  access to the law of the hitting time of a square-root boundary by a Bessel process. This extends classical results by Breiman and Shepp, concerning   Brownian motion, and recovers by different means, extensions for Bessel processes, obtained independently by Delong and Yor.
\bigskip
\bigskip


 Let $B_t$ be the standard real valued Brownian motion and for $\nu>0$, introduce the geometric Brownian motion $\E_t^{(-\nu)}$ and its exponential functional $\A_t^{(-\nu)}$
 $$\E_t^{(-\nu)}:= exp (B_t-\nu t)$$
 $$\A_t^{(-\nu)}:=\int_0^t (\E_s^{(-\nu)})^2 ds.$$
 
 Lamperti's representation theorem \cite{lamperti} applied to $\E_t^{(-\nu)}$ states
\begin{equation}\label{lamperti}
\E_t^{(-\nu)}=R_{\A_t^{(-\nu)}}^{(-\nu)}
\end{equation}
 where $(R_u^{(-\nu)}, \, u\leq T_0(R^{(-\nu)}))$ denotes the Bessel process of index $(-\nu)$ (equivalently of dimension $\delta=2(1-\nu)$), starting at 1, which is an $\R_+$-valued diffusion with infinitesimal generator $\L^{(-\nu)}$ given by
 $$\L^{(-\nu)}f(x)={1\over2} f''(x) + {1-2\nu\over 2x} f'(x), \quad f\in C_b^2(\R_+^\star).$$


Let us  remark that, in the special case $\nu=1/2$, equation (\ref{lamperti}) is nothing else but the Dubins-Schwarz representation of the exponential martingale $\E_t^{(-1/2)}$ as Brownian motion time changed with $\A_t^{(-1/2)}$.

For a short summary of relations between Bessel processes and exponentials of Brownian motion, see e.g.  Yor \cite{yor1}.

Let us consider now the following random variable $Z$, which is often called a perpetuity in the mathematical finance literature:
$$Z:=\A_\infty^{(-\nu)}=\int_0^\infty (\E_s^{(-\nu)})^2 ds$$
We deduce directly from  (\ref{lamperti}) that
$$\A_\infty^{(-\nu)}=T_0(R^{(-\nu)})$$
where $T_0:=\inf\{u : X_u=0\}$, and it is well-known
(see \cite{dufresne}, \cite{yor2}),  that 
\begin{equation}\label{dufresne}
\A_\infty^{(-\nu)}\, {\stackrel{(law)} {=}}\,{1\over 2\gamma_\nu}
\end{equation}
where $\gamma_\nu$ is a gamma variable with parameter $\nu$ (i.e. with density ${1\over \Gamma(\nu)}x^{\nu-1}e^{-x}{\bf 1}_{\R_+}$).

Our  main result characterizes the law of the hitting time of a parabolic boundary by $R_u^{(-\nu)}$ which corresponds to a Bessel process of dimension $d<2$.
\begin{theorem}
Let  $0<b<c$, and $\sigma:=\inf\{ u : (R_u^{(-\nu)})^2={1\over c}(b+u)\}$ with $R_0^{(-\nu)}=1$. 
\begin{equation}
\label{theorem1}
E[(b+\sigma)^{-s}]=c^{-s} {E[(1+2b\gamma_{\nu+s})^{-s}]\over E[(1+2c\gamma_{\nu+s})^{-s}]}, \quad \hbox{\it for any}\,\,\, s\geq0
\end{equation}
\end{theorem}

Proof: 
using the strong Markov property and the stationarity of the increments of  Brownian motion, we obtain that for any stopping time $\tau$ of the Brownian motion 
$$\A_\infty^{(-\nu)}=:Z=\A_\tau ^{(-\nu)}+(\E_\tau^{(-\nu)} )^2Z' $$
where $Z'$ is independent of $(\A_\tau^{(-\nu)}, \E_\tau^{(-\nu)})$ and $Z\,{\stackrel{(law)} {=}}\,Z'$.

This implies,  by (\ref{lamperti}), that $Z$ satisfies the following affine equation (see \cite{vervaat} for a survey about these equations)
\begin{equation}
\label{renewal1}
\A_\infty^{(-\nu)}=:Z= \A_\tau ^{(-\nu)}+(R^{(-\nu)}_{\A_\tau^{(-\nu)}})^2 Z'
\end{equation} 
where $Z'$ is independent of $(\A_\tau^{(-\nu)}, R^{(-\nu)}_{\A_\tau^{(-\nu)}})$ and $Z\,{\stackrel{(law)} {=}}\,Z'$.

Obviously, $\sigma<T_0(R^{(-\nu)})$. 
Taking now :
$$ \tau =\inf\{ t : (R^{(-\nu)}_{\A_t^{(-\nu)}})^2 ={1\over c}(b+ \A_t^{(-\nu)}) \}$$
we get   $\A_\tau^{(-\nu)}=\sigma$, and the identity in law 
\begin{equation}
\label{renewal2}
b+Z\,{\stackrel{(law)} {=}}\,(b+\sigma)(1+{Z\over c})\end{equation} 
where the variables $\sigma$ and $Z$ on the right-hand side are independent.

As a result, we obtain the Mellin-Stieltjes transform of $\sigma$:
$$E[(b+\sigma)^{-s}]=c^{-s} {E[(b+Z)^{-s}]\over E[(c+Z)^{-s}]}$$
But, from (\ref{dufresne}) 
$$E[(b+\sigma)^{-s}]=c^{-s} {E[(2\gamma_\nu)^s{1\over (1+2b\gamma_\nu)^s}]\over 
E[(2\gamma_\nu)^s{1\over (1+2c\gamma_\nu)^s}]}$$
which gives the result.\qed

One can now use the duality between the laws of Bessel processes of dimension $d$ and $4-d$ to get the analogous result of Theorem 1, and recover the result of Delong \cite{delong1}, 
\cite{delong2}, and Yor \cite{yor} which deals with  the case  $d\geq2$.

\begin{theorem}
Let  $0<b<c$, and $\sigma:=\inf\{ u : (R_u^{(\nu)})^2={1\over c}(b+u)\}$ with $R_0^{(\nu)}=1$. 
\begin{equation}
\label{theorem2}
E[(b+\sigma)^{-s}]=c^{-s} {E[(1+2b\gamma_{s})^{-s+\nu}]\over E[(1+2c\gamma_{s})^{-s+\nu}]}, \quad \hbox{\it for any}\,\,\, s\geq0.
\end{equation}

\end{theorem}

Proof : it is based on the following classical relation between the laws of the Bessel processes with indices $\nu$ and $-\nu$: 
\begin{equation}
\label{doob}
{{\mathcal P}_x^{(\nu)}}_{\vert{\mathcal F}_t}={(X_{t\wedge T_0})^{2\nu}\over x^{2\nu}} . 
{{\mathcal P}_x^{(-\nu)}}_{\vert{\mathcal F}_t}
\end{equation}
which implies that 
$$E_1^{(\nu)}[(b+\sigma)^{-s}]=E_1^{(-\nu)}[X_\sigma^{2\nu}(b+\sigma)^{-s}]={1\over c^\nu}E_1^{(-\nu)}[(b+\sigma)^{-s+\nu}]$$
Theorem 1 gives the result. \qed

Finally, it is easily shown, thanks to the classical representations of the Whittaker functions (see Lebedev \cite{lebedev}), that the right-hand sides of (\ref{theorem1}) and (\ref{theorem2}) are expressed in terms of ratios of Whittaker functions.

{\bf Aknowledgement:} We would like to thank Daniel Dufresne for useful and enjoyable  discussions on the subject.

\end{document}